\documentclass[12pt,oneside,english]{amsart}
\usepackage[letterpaper]{geometry}
\usepackage{babel}
\usepackage{amssymb}
\usepackage{mathrsfs}
\usepackage[arrow,matrix]{xy}
\usepackage[symbol]{footmisc}
\usepackage{mathptmx}
 \theoremstyle{plain}    
 \newtheorem{thm}{Theorem} %[section]
\theoremstyle{definition}
 \newtheorem{defn}[thm]{Definition}
 \newtheorem{claim}{Claim}
\theoremstyle{plain}    
 \newtheorem{prop}[thm]{Proposition} %%Delete [thm] to re-start numbering
 \newtheorem{cor}[thm]{Corollary} %%Delete [thm] to re-start numbering
 \theoremstyle{definition}
\theoremstyle{remark}

\def\ie{i.e.\ }

\def\N{\mathbb{N}}
\def\Z{\mathbb{Z}}

\def\C{\mathbb{C}}

\def\B{\mathscr{B}}

\def\M{\mathcal{M}}

\DeclareMathOperator{\spec}{spec}
\DeclareMathOperator{\image}{image}

\DeclareMathOperator{\rank}{rank}

\let\tens\otimes

\def\tmin{\otimes_{\rm{min}}}

\def\incl{\mathrel{\hookrightarrow}}

\def\L{\mathscr{L}}
\let\amalg\sqcup

\def\UC{C^*_u}

\DeclareMathOperator{\Alt}{Alt}

\title{Non-$K$-exact uniform Roe C*-algebras}

\author{J\'an \v Spakula}
\address{Mathematisches Institut, Universit\"at M\"unster, Einsteinstr.\ 62,
    48149 M\"unster, Germany}
\email{jan.spakula@uni-muenster.de}

\begin{document}

\begin{abstract}
  We prove that uniform Roe C*-algebras $C^*_uX$ associated to some expander
  graphs $X$ coming from discrete groups with property ($\tau$) are not
  $K$-exact. In particular, we show that this is the case for the expander
  obtained as Cayley graphs of a sequence of alternating groups
  (with appropriately chosen generating sets).\\[1ex]


\end{abstract}

\maketitle

\section{Introduction}

Uniform Roe C*-algebras (also called uniform translation C*-algebras)
provide, among other things, a link between coarse geometry and C*-algebra
theory via the following theorem, which connects a coarse--geometric property
of a discrete group $\Gamma$ with a purely analytic property of its reduced
C*-algebra $C^*_r\Gamma$:

\begin{thm}[Guentner--Kaminker \cite{guent-kam:exactness} and Ozawa
  \cite{ozawa:exactness}] Let $\Gamma$ be a finitely generated discrete
  group. Then the following are equivalent:
  \begin{itemize}
    \item $\Gamma$ has property A (see \cite{yu:CBC-for-embeddable}),
    \item the reduced group C*-algebra $C^*_r\Gamma$ is exact,
    \item the uniform Roe C*-algebra $\UC|\Gamma|$ is nuclear.
  \end{itemize}
\end{thm}

The only known examples of groups which do not have property A are Gromov's
random groups \cite{gromov:spacesandquestions}. The theorem also
characterizes nuclearity of uniform Roe C*-algebras of discrete
groups.

More generally, one may ask if an analogue of the above theorem for general
bounded geometry metric spaces $X$ is true, i.e.\ if property A of
$X$ is equivalent to nuclearity of $\UC X$. One proof was obtained by
representing $\UC X$ as a groupoid C*-algebra (see \cite{ska-tu-yu}) and
referring to the general groupoid C*-algebra theory (see
\cite{ana-ren:groupoids}). For a more elementary argument for one direction
see \cite{roe:LOCG}.

In further attempt to generalize, we may ask about a $K$-theoretic analogue
of C*-algebraic property of nuclearity. Namely, we may seek some (coarse)
geometric conditions of $X$ that would imply $K$-nuclearity of uniform
Roe C*-algebras. The first step beyond the realm of property A is
coarse embeddability into a Hilbert space. Using groupoid language, Ulgen
proved that if $X$ admits a coarse embedding into a Hilbert space, then $\UC
X$ is $K$-nuclear (see \cite[proof of Theorem 3.0.12]{ulgen:thesis}).  The
argument uses the fact that groupoids with the Haagerup property are
$K$-amenable, and that crossed products of $K$-nuclear algebras with
$K$-amenable groupoids are again $K$-nuclear (see \cite{tu:BCamenable} and
\cite{tu:BCgroupoids}).

\smallskip

On the other side, we may look for examples of spaces, whose uniform
Roe C*-algebras are not nuclear. An unpublished result of Higson
asserts that uniform Roe C*-algebras of expander graphs constructed
from groups with property (T) are not exact, and therefore not nuclear.

Wandering into the $K$-theoretic territory, we can ask for examples of
spaces, whose uniform Roe C*-algebras are not $K$-nuclear.  Ulgen
defined $K$-exactness in \cite{ulgen:thesis} as a generalization of exactness
in the context of $K$-theory. She also proved that if a separable C*-algebra
is $K$-nuclear, then it is also $K$-exact. Unfortunately, this cannot be
applied to uniform Roe C*-algebras, which are usually not separable.
In this paper, we show that for some spaces $X$, $\UC X$ is not even
$K$-exact. Our examples are expander graphs, constructed out of groups with
property $(\tau)$ with respect to a family of subgroups $\L$ (see
\cite{lubotzky:prop-tau}). Certain assumption on the family $\L$ is
required at this point:

\begin{thm}\label{thm:main}
Let $\Gamma$ be a finitely generated discrete group with property
$\tau(\L)$. Assume that
\begin{itemize}
  \item[($\star$)] $\Gamma$ has $\tau(\L')$, where $\L'=\{N_1\cap N_2\mid
    N_1,N_2\in \L\}$.
\end{itemize}
Then $\UC X$ is not $K$-exact, where $X=\amalg_{N\in \L}\Gamma/N$.
\end{thm}

Using results of Kassabov \cite{kassabov:sym-groups}, we obtain the following
corollary:

\begin{cor}\label{cor:alt}
  There is a sequence $(n_i)_{i\in\N}$, such that the uniform Roe C*-algebra
  of the expander obtained as a coarse disjoint union of Cayley graphs of the
  alternating groups $\Alt(n_i)$ (with appropriately chosen generating sets)
  is not $K$-exact.
\end{cor}

The question of $K$-exactness for uniform Roe C*-algebras of
expander graphs is closely related to the same question for C*-algebras of
the type $\prod_q\M_q(\C)$. This has been settled negatively in various
contexts by Ozawa \cite{ozawa:SL3Z} and by Manuilov--Thomsen
\cite{manuilov-thomsen:ext}. Both constructions extend the work of Wassermann
\cite{wassermann:T}. 

Projections in (uniform) Roe C*-algebras similar to the one that is used in
the construction in this paper are called \emph{ghosts} and were studied in
the context of the ideal structure of Roe C*-algebras 
\cite{chen-wang:ideal-structure,qin:ghosts-in-perturbed-expanders}.

\smallskip

The paper is organized as follows: In the next section, we recall the
definitions of properties and objects involved. Section \ref{sec:thm} is
devoted to the proof of Theorem \ref{thm:main} and in the last section we
show how to deduce Corollary \ref{cor:alt}.

\textbf{Acknowledgment:} The author would like to thank Guoliang Yu for
helpful and enlightening conversations and never-ending encouragement and
support, the reviewer of the previous version of this paper for enormous
simplification and generalization of previous argument, and finally Martin
Kassabov for helpful comments about the alternating groups.

\section{Definitions}

\subsection{$K$-exactness}
Recall that a C*-algebra $A$ is exact, if $\cdot\tmin A$ is an exact functor,
\ie if we $\min$--tensor every term in a short exact sequence with $A$, the
sequence stays exact. If it does, then we obtain a 6-term exact sequence in
$K$-theory (as below). It may happen that the tensored short exact sequence
is not exact for a non--exact C*-algebra, but there is always an exact 6-term
sequence in $K$-theory. Let us be more precise:
\begin{prop}[{\cite[2.3.2]{ulgen:thesis}}]
  For a C*-algebra $A$, the following are equivalent:
  \begin{itemize}
  \item for any exact sequence of C*-algebras $0\to I\to B\to B/I\to 0$,
    there is a cyclic 6-term exact sequence in $K$-theory:
    $$
    \xymatrix{
      K_0(I\tmin A)\ar[r] & K_0(B\tmin A)\ar[r] & K_0(B/I\tmin A)\ar[d] \\
      K_1(B/I\tmin A)\ar[u] & K_1(B\tmin A)\ar[l] & K_1(I\tmin A).\ar[l] }
    $$
    (Note that in general, there might be no such 6-term $K$-theory sequence
    at all.),
  \item for any exact sequence of C*-algebras $0\to I\to B\to B/I\to 0$, the
    sequences
    $$
    K_i(I\tmin A)\to K_i(B\tmin A)\to K_i(B/I\tmin A),
    $$
    are exact in the middle for both $i=0,1$.
  \end{itemize}
\end{prop}
\begin{defn}
  We say that a C*-algebra $A$ is \emph{$K$-exact}, if it satisfies the
  conditions in the previous proposition.
\end{defn}
For separable C*-algebras, a sufficient condition for $K$-exactness is
$K$-nuclearity. The argument \cite[proposition 3.4.2]{ulgen:thesis} can be
summarized as follows: the $\max$--tensor product always preserves exact
sequences, and if a C*-algebra is $K$-nuclear, then $\min$--tensor products
and $\max$--tensor products with it are $KK$--equivalent. However, this
relies on the key properties of $KK$-theory (existence and associativity of
Kasparov product), which have been proved only for separable C*-algebras in
general.

\subsection{Uniform Roe C*-algebras}

A metric space $X$ has bounded geometry, if for each $r>0$ the number of
points in any ball or radius $r$ is uniformly bounded. We say that $X$ is
uniformly discrete, if there exists $c>0$, such that any two distinct points
of $X$ are at least $c$ apart.
\begin{defn}
Let $X$ be a uniformly discrete metric space with bounded geometry.
We say that an $X$-by-$X$ matrix $(t_{yx})_{x,y\in X}$ with complex entries
has \emph{finite propagation}, if there exists $R\geq0$, such that $t_{yx}=0$
whenever $d(x,y)\geq R$. We say that such a matrix is \emph{uniformly
bounded}, if there exists $T\geq0$, such that $|t_{yx}|\leq T$ for all $x,y\in
X$.
\end{defn}
Let $A(X)$ be the algebra of all finite propagation matrices which are
uniformly bounded. It is easy to see that each element of $A(X)$ represents a
bounded operator on $\ell^2(X)$ (see \cite[lemma 4.27]{roe:LOCG}). This
yields a representation $\lambda:A(X)\to\B(\ell^2(X))$.
\begin{defn}
The \emph{uniform Roe C*-algebra} $\UC X$ of $X$ is defined to be
the norm closure of $\lambda(A(X))\subset \B(\ell^2(X))$.
\end{defn}

\subsection{Expanders and property $(\tau)$}
\begin{defn}
  An \emph{expander} is a sequence $X_n$ of finite graphs with the
  properties:
  \begin{itemize}
    \item The maximum number of edges emanating from any vertex is uniformly
      bounded.
    \item The number of vertices of $X_n$ tends to infinity as $n$ increases.
    \item The first nonzero eigenvalue of the Laplacian, $\lambda_1(X_n)$, is
      uniformly bounded away from zero, say by $\lambda>0$.
  \end{itemize}
  We think of $X_n$ as of a discrete metric space, where the points are
  vertices of the graph, and the metric is given by the path distance in the
  graph. We understand the sequence as one metric space $\amalg_nX_n$ via
  the coarse disjoint union construction.
\end{defn}
Let us recall one possibility of how to construct a coarse disjoint union of
finite spaces: Given a sequence $(X_q)_{q\in\N}$ of finite metric spaces, we
define their coarse disjoint union $\amalg_qX_q$ to be the set $\cup_qX_q$
endowed with the metric inherited from individual $X_q$'s together with the
condition $d(X_q,X_{q'})={\max(q,q')}$ for $q\not=q'$.

The first explicit examples of expanders were constructed by Margulis as
$\amalg_q\Gamma/\Gamma_q$, where $\Gamma$ is a finitely generated group with
property (T) (with a fixed generating set), and $\Gamma_q\leq \Gamma$ is a
decreasing sequence of normal subgroups with finite index, such that
$\bigcap_q\Gamma_q=\{1\}$. This construction eventually led to Lubotzky's
property ($\tau$) \cite{lubotzky:prop-tau}.

\begin{defn}
  Let $\Gamma$ be a finitely generated group and $\L$ a countable family of
  finite index normal subgroups of $\Gamma$. We also assume that $\L$ is
  infinite, and that $[\Gamma:N]\to\infty$ as $N\to\infty$\footnote[2]{By
    $N\to\infty$ we mean ``outside the finite sets''.}, $N\in\L$. We say that
  $\Gamma$ has \emph{property ($\tau$) with respect to the family $\L$}
  (written also $\tau(\L)$) if the trivial representation is isolated in the
  set of all unitary representations of $\Gamma$, which factor through
  $\Gamma/N$, $N\in\L$. We say that $\Gamma$ has \emph{property ($\tau$)} if
  it has this property with respect to the family of all finite index normal
  subgroups.
\end{defn}

Property $\tau(\L)$ is equivalent to $\amalg_{N\in\L}\Gamma/N$ being an
expander \cite[Theorem 4.3.2]{lubotzky:prop-tau}.

\section{Proof of Theorem \ref{thm:main}}
\label{sec:thm}

We describe the construction starting with $\Gamma$, a finitely generated
discrete group and a countable family $\L$ of normal subgroups of $\Gamma$
with finite index. We also fix a finite symmetric generating set $S$ of
$\Gamma$.

For each $N\in \L$, we denote $G_N=\Gamma/N$ and by $q_N:\Gamma\to G_N$ the
quotient map.  We let $X$ to be a coarse disjoint union of the Cayley graphs
of $G_N$'s with respect to the generating sets $\{ q_N(g)\mid g\in S\}$.

Let $\lambda_N:G_N\to \B(\ell^2G_N)$ be the left regular representation of
$G_N$.  Denote also $\tilde\lambda_N=\lambda_N\circ q_N:\Gamma \to
\B(\ell^2G_N)$ and
$\Lambda=\oplus_{N\in\L}\tilde\lambda_N:\Gamma\to\B(\ell^2X)$.

\begin{claim}
  For each $N\in\L$, we can choose an irreducible representation $\pi_N:G_N\to
  \B(H_N)$, so that
  \begin{itemize}
  \item[($\star\star$)] $\dim(H_N)\to \infty$ as $N\to \infty$.
  \end{itemize}
\end{claim}
\begin{proof}
  We shall use the fact that if $\Gamma$ has $\tau(\L)$, then for each fixed
  $d>0$, there are only finitely many non-equivalent irreducible
  $d'$-dimensional ($d'\leq d$) representations of $\Gamma$ factoring through
  some $G_N$, $N\in\L$. The same conclusion is known for groups with property
  (T) \cite{wang:on-isolated-points75,
    delaharpe-robertson-valette:on-the-spectrum93}, the argument for groups
  with $(\tau)$ is outlined also in \cite[Theorem 3.11]{lubotzky-zuk:tau}.
  
  Now assume that the claim doesn't hold, that is, there is a sequence
  $N_n\in\L$, such that all irreducible representations of $G_{N_n}$'s are at
  most $d$-dimensional. By the above fact, all of them come from a finite set
  of irreducible representations $\{\rho_1,\dots,\rho_m\}$ of $\Gamma$.
  Consequently, each $G_{N_n}$ embeds as a subgroup of
  $K=\image(\rho_1\oplus\cdots\oplus\rho_m)$. Since every $\rho_i$ factors
  through some $G_{N_k}$, its image is a finite group. Hence $K$ is
  finite, so we obtain a contradiction with the assumption that $|G_{N_n}|\to
  \infty$ as $n\to\infty$.
\end{proof}

Denote $\tilde
\pi_N=\pi_N\circ q_N$, $H=\oplus_{N\in\L}H_N$ and
$\pi=\oplus_{N\in\L}\pi_N\circ q_N:\Gamma\to\B(H)$. Let us summarize the
notation in the following diagram:
$$
\xymatrix{ \Gamma \ar[r]^{q_N} \ar[rd]_{\tilde\lambda_N} &
  G_N\ar[d]^{\lambda_N} \\ & \B(\ell^2G_N)}
\qquad\quad
\xymatrix{ \Gamma \ar[r]^{q_N} \ar[rd]_{\tilde\pi_N} &
  G_N\ar[d]^{\pi_N} \\ & \B(H_N)}
$$

Finally, we let $B=\prod_{N\in\L}\B(H_N)$ and $J=\oplus_{N\in\L}\B(H_N)$. We
obtain an exact sequence of C*-algebras
$$
0\to J\to B\to B/J\to 0.
$$
We shall use this sequence to show that $\UC X$ is not $K$-exact. We
construct a projection $e\in \UC X\tens B$, whose $K_0$-class will violate
the exactness of the $K$-theory sequence
$$
K_0(\UC X\tens J)\to K_0(\UC X\tens B)\to K_0(\UC X\tens (B/J)).
$$
To construct such $e$, we let
$$
T=\frac1{|S|}\sum_{g\in S}(\Lambda\tens \pi)(g) \in \UC X \tens B\subset
\B(\ell^2X\tens H).
$$
A few remarks are in order:
\begin{itemize}
\item If we denote $s=\frac1{|S|}\sum_{g\in S}g\in\C\Gamma$, then
  $T=(\Lambda\tens\pi)(s)$.
\item $T$ is ``diagonal'' with respect to the decomposition $\ell^2X\tens H
  = \oplus_{M,N\in\L}\ell^2G_N\tens H_M$. This is clear from the fact that
  each $\ell^2G_N\tens H_M$ is a $\Lambda\tens\pi$--invariant subspace. We
  denote its ``entries'' by $T_{NM}\in \B(\ell^2G_N\tens H_M)$.
\item With our choice of metric on $X$, $\Lambda(g)\in \B(\ell^2X)$ has
  propagation $1$ for each $g\in S$.
\end{itemize}
The construction is finished by proving three claims, which we state and give
some remarks about them. The proofs are spelled out afterward.
\begin{claim}\label{claim:1-is-isolated}
  $1\in \spec(T)$ is an isolated point.
\end{claim}
For proving this, we need to use some form of property $\tau$. The necessary
condition for Claim \ref{claim:1-is-isolated} is that $\Gamma$ has
$\tau(\L)$. However, this by itself is not sufficient, since the Claim
requires uniform bound on the spectral gap for all $\lambda_N\tens\pi_M$, not
just $\lambda_N$'s. The condition $(\star)$ ensures this.

Claim \ref{claim:1-is-isolated} allows us to define the projection $e\in \UC
X\tens B$ to be the spectral projection of $T$ corresponding to $1\in
\spec(T)$. Note that $e$ is also ``diagonal'' as $T$, hence we can decompose
it into projections $e_{NM}\in \B(\ell^2G_N\tens H_M)$. It is clear from the
definition of $T$ that each $e_{NM}$ is in fact the projection onto the
subspace of $\Gamma$-invariant vectors in $\ell^2G_N\tens H_M$. In fact, if
$\Gamma$ has property (T), then $e$ is the image of the Kazhdan projection
$p_0\in C^*_{\rm max}\Gamma$ under $\Lambda\tens\pi$.
\begin{claim}\label{claim:goes-to-0}
  $e$ maps to $0\in \UC X\tens (B/J)$.
\end{claim}
The key observation here is that if for a fixed $N$, $e_{NM}$ are eventually
$0$, then Claim \ref{claim:goes-to-0} holds. This is where the condition
($\star\star$) is used.

\begin{claim}\label{claim:doesnt-come-from-J}
  $[e]\in K_0(\UC X\tens B)$ does not come from any class in $K_0(\UC X\tens
  J)$.
\end{claim}
This claim is proved by ``detecting the diagonal'' of $e$. More precisely,
observe that $e_{NN}\not=0$ for every $N\in\L$, since $\pi_N$ is conjugate to
a subrepresentation of $\lambda_N$ and hence $\lambda_N\tens\pi_N$ has
nonzero invariant vectors. However, any element coming from $\UC X\tens J$
will have the $NN$-entries eventually $0$. The construction of a
*-homomorphism that detects this is essentially due to Higson
\cite{hig-laf-ska:counters}.

\begin{proof}[Proof of Claim \ref{claim:1-is-isolated}]
  Taking $N\in\L$, the $G_N$-action on $\ell^2G_N\tens H_N$ is via
  $\lambda_N\tens \pi_N$. This representation contains the trivial
  representation (since $\lambda_N$ contains the conjugate of $\pi_N$, as it
  does any irreducible representation of $G_N$), so there are nonzero
  $G_N$-invariant vectors in $\ell^2G_N\tens H_N$. Therefore, $1\in
  \spec(T_{NN})$.

  To show that $1$ is actually isolated in each $\spec(T_{NM})$ with the
  uniform bound on the size of the gap, we shall use the condition
  ($\star$). Property $\tau(\L')$ says that we have such a uniform bound on
  the size of the spectral gap of the image of $\rho(T)$ for all the
  representations $\rho$ of $\Gamma$ which factor through some of $G_L$,
  $L\in \L'$ \cite[Theorem 4.3.2]{lubotzky:prop-tau}. Using that
  $\tilde\pi_M$ is contained in $\tilde\lambda_M$, we obtain
  $$
    \ker(\tilde\lambda_N\tens\tilde\pi_M)=\ker(\tilde\lambda_N)\cap
    \ker(\tilde \pi_M)\supseteq
    \ker(\tilde\lambda_N)\cap\ker(\tilde\lambda_M) = N\cap M.
  $$
  This shows that $\tilde\lambda_N\tens\tilde\pi_M$ factors through
  $\Gamma/(N\cap M)$ and the proof is finished.
\end{proof}

\begin{proof}[Proof of Claim \ref{claim:goes-to-0}]
  Denote $A=\prod_{N\in\L}\B(\ell^2G_N)$, a product of matrix algebras.  It
  is clear that $T\in (\UC X\tens B)\cap(A\tens B)$, and so also $e\in A\tens
  B\subset \B(\ell^2X\tens H)$. 

  For $N\in\L$, let us examine the $\B(\ell^2G_N)\tens\B(H)$--component of
  $e$. Denote by $P_N\in\B(\ell^2X)$ the projection onto $\ell^2G_N$. It
  suffices to show that $e_N=(P_N\tens 1_B)e(P_N\tens 1_B)\in
  \B(\ell^2G_N\tens H)$ actually belongs to $\B(\ell^2G_N)\tens J$, since
  that shows that $e\in \UC X\tens J$, and therefore maps to $0\in \UC
  X\tens(B/J)$.

  Further decompose $e_N$ into $e_{NM}\in \B(\ell^2G_N\tens H_M)$. Recall
  that $e_{NM}\not=0$ if and only if $\ell^2G_N\tens H_M$ has nonzero invariant
  vectors. Since the representation $\tilde\pi_M$ is irreducible, this is
  further equivalent to $\tilde\pi_M$ being conjugate to a subrepresentation
  of $\tilde \lambda_N$. But by ($\star\star$), this can only happen for
  finitely many $M$'s, since $\tilde\lambda_N$ is fixed and
  $\dim(H_M)\to\infty$.
\end{proof}

\begin{proof}[Proof of Claim \ref{claim:doesnt-come-from-J}]
  For $N\in\L$, denote $C_N=\B(\ell^2G_N\tens H_N)$. We construct a
  *-ho\-mo\-mor\-phism $f:\UC X\tens B\to \prod_NC_N\big/\oplus_NC_N$, such that
  $f_*([e])\not=0\in K_0(\prod_NC_N\big/\oplus_NC_N)$, but $f_*([x])=0$ for any
  $[x]\in K_0(\UC X\tens J)$.

  We first embed $\UC X$ into a direct limit of C*-algebras $A_k$, defined
  below. We enumerate $\L=\{N_k\mid k\in\N\}$ and put
  $A^q=\B(\ell^2G_{N_q})$, $A^{0k}=\B(\ell^2(\amalg_{q\leq k}G_{N_q}))$ and
  finally $A_k=A^{0k}\oplus \prod_{q>k}A^q$, $k\geq 1$. There are obvious
  inclusion maps $A_k\incl A_l$ for $k<l$, so we can form a direct limit
  $A_0=\lim_kA_k$.  It follows from the condition on distances
  $d(G_{N_q},G_{N_p})$ that each finite propagation operator on $\ell^2X$ is
  a member of some $A_k\subset\B(\ell^2X)$, hence $\UC X\incl A_0$.

  For $k\geq1$, denote $B^k=\B(H_{N_k})$ (so that $B=\prod_{k\in\N}B^k$) and
  define $f_{k}$ as the following composition:
  \begin{align*}
  A_k\tens B &= \Bigl(A^{0k}\oplus \prod_{q>k}A^q\Bigr)\tens B\incl
  \left(A^{0k}\tens B\right)\oplus
    \prod_{q>k}\left(\left[A^q\tens B^q\right]\oplus
      \biggl[A^q\tens\Bigl(\prod_{p\not=q}B^p\Bigr)\biggr]\right)
  \stackrel{\mathrm{proj}}{\longrightarrow}\\
  &\stackrel{\mathrm{proj}}{\longrightarrow} \prod_{q>k}A^q\tens B^q = 
  \prod_{q>k}C_{N_k}
  \stackrel{\mathrm{quot}}{\longrightarrow} \prod_qC_{N_q}\Big/\sum_qC_{N_q}.
  \end{align*}
  It is easy to see that $f_k$'s commute with inclusions
  $A_k\tens B \incl A_l\tens B$, $k<l$.
  Consequently, we obtain a *-homomorphism $f:\UC X\tens(B\rtimes\Gamma)\to
  \prod_N C_N\big/\oplus_N C_N$.
  
  It is known that $K_0(\prod_N C_N\big/\oplus_N C_N)$ embeds into
  $\prod_N\Z\big/\oplus_N\Z$. Examining the construction of $f$, we see that
  $f_*([e])$ is the class of the sequence $k\mapsto\rank(e_{N_kN_k})$ in
  $\prod_{N_k}\Z\big/\oplus_{N_k}\Z$. As already noted, every term of this
  sequence is nonzero.

  On the other hand, any projection $p\in \UC X\tens J$ has only finitely
  many nonzero $C_N$-components, and so $f_*([p])=0$ in
  $\prod_N\Z\big/\oplus_N\Z$. This obviously extends to the whole of $K_0(\UC
  X\tens J)$.
\end{proof}

\section{Expanders coming from finite groups}

Start with a sequence $G_n=\langle S_n\rangle$, $n\in\N$, of finite groups
with chosen generating sets $S_n$ of some fixed size. Let $\Gamma$ be a
subgroup of $G=\prod_{n\in\N}G_n$ generated by a finite set $S\subset G$ that
projects onto $S_n$ in each factor. Denote by $q_n:\Gamma\to G_n$ the natural
projection, $N_n=\ker(q_n^{-1})\subset \Gamma$ and $\L=\{N_n\mid n\in\N\}$.

Now $\Gamma$ has $\tau(\L)$ if and only if the Cayley graphs of $G_n$'s with
respect to $S_n$'s constitute an expander. In order to apply Theorem
\ref{thm:main}, we need to verify the condition ($\star$). It seems to be
open whether $\tau(\L)$ implies ($\star$) in general. For instance, it would
be sufficient to know a positive answer to \cite[Question
1.14]{lubotzky-zuk:tau}, that is, whether $\tau(\L)$ implies $\tau(\L'')$,
where $\L''$ is the closure of $\L$ under finite intersections. See also
\cite[Question 6]{kassabov-nikolov:cartesian-products06} for a discussion
when there is a finitely generated dense $\Gamma\subset\prod_{n\in\N}G_n$
which has property ($\tau$) (with respect to all finite index subgroups).

To prove Corollary \ref{cor:alt}, we appeal to the result of Kassabov
\cite[p.\ 352]{kassabov:sym-groups}, which says that there is a finitely
generated dense subgroup $\Gamma\subset \prod_s\Alt\left(
  (2^{3s}-1)^6\right)$ which has property ($\tau$). Hence taking a finite
symmetric generating set of $\Gamma$ projects to generating sets of the
individual factors, making their Cayley graphs into an expander. Moreover,
the condition ($\star$) is obviously satisfied. Consequently, we have shown
Corollary \ref{cor:alt}.

%\bibliography{biblio.bib}
%\bibliographystyle{siam}
%\printbibliography

\end{document}